\documentclass[12pt]{amsart}
\usepackage{amssymb,latexsym, amsthm,epsf}

\newtheorem{thm}{Theorem}[section] 
\newtheorem{defn}[thm]{Definition}

\newtheorem{cor}[thm]{Corollary}

\newtheorem{prob}[thm]{Problem}

\newtheorem{prop}[thm]{Proposition}
\newtheorem{rem}[thm]{Remark}

\def\ie{{i.e.}} 

\def\C{{\mathbb C}}
\def\F{{\mathbb F}}

\def\R{{\mathbb R}}
\def\P{{\mathbb P}}
\def\PP{{\mathbb P}}

\def\Z{{\mathbb Z}}

\def\ra{{\rightarrow}}

\def\({\left(}
\def\){\right)}

\newcommand{\set}[1]{{\{#1\}}}
\newcommand{\power}[1]{{\left|{#1}\right|}}

\newcommand\suchthat{{\,:\,}}

\newcommand\defin[1]{{\it{#1}}}
\long\def\forget#1\forgotten{}

\newcommand\Lpair[2]{{\left[{#1},{#2}\right]}}
\newcommand{\sbs}{\subset}

\newcommand{\Int}{\mathop {\rm Int}}
\newcommand{\cB}{{\mathcal B}}
\newcommand{\be}{\beta}

\newcommand{\LPs}{Lefschetz pairs}
\def\fg{fundamental group}
\renewcommand\th[1]{{$#1$th}}
\newcommand{\ba}{{\bf a}}

\newcommand{\G}{\Gamma}
\newcommand{\Ga}{\Gamma}

\newif \iffurther 
\furtherfalse

\newif \iffurther 
\furtherfalse

\newif\ifXY 
\XYfalse    

\ifXY
\usepackage{xy}
\fi
\ifXY
\xyoption{all}
\fi

\begin{document}

\title[Conjugation-free presentations of arrangements]{A conjugation-free geometric presentation of fundamental groups of arrangements}

\author[Eliyahu, Garber, Teicher]{Meital Eliyahu$^1$, David Garber$^1$ and Mina Teicher}

\stepcounter{footnote}
\footnotetext{Partially supported by a
grant from the Ministry of Science, Culture and Sport, Israel and
the Russian Foundation for Basic research, the Russian
Federation.}

\address{Meital Eliyahu, Department of Mathematics, Bar-Ilan University, 52900 Ramat-Gan, Israel}
\email{eliyahm@macs.biu.ac.il}

\address{David Garber, Department of Applied Mathematics, Faculty of
  Sciences, Holon Institute of Technology, 52 Golomb st., PO
  Box 305, 58102 Holon, Israel}
\email{garber@hit.ac.il}

\address{Mina Teicher, Department of Mathematics, Bar-Ilan University, 52900 Ramat-Gan, Israel}
\email{teicher@macs.biu.ac.il}
\keywords{}

\begin{abstract}
We introduce the notion of a {\em conjugation-free geometric presentation} for a fundamental group of a line arrangement's complement, and we show that the fundamental groups of the following family of arrangements have a conjugation-free geometric presentation: A real arrangement $\mathcal L$, whose graph of multiple points is a union of disjoint cycles, has no line with more than two multiple points, and where the multiplicities of the multiple points are arbitrary.

We also compute the exact group structure (by means of a semi-direct product of groups) of the arrangement of $6$ lines whose graph consists of a cycle of length $3$, and all the multiple points have multiplicity $3$.
\end{abstract}

\maketitle

\section{Introduction}
The fundamental group of the complement of plane curves is a very important topological invariant, which can be also computed for line arrangements. We list here some applications of this invariant.

Chisini \cite{chisini}, Kulikov \cite{Kul,Kul2} and Kulikov-Teicher \cite{KuTe} have used the fundamental group of complements of branch curves of
generic projections in order to distinguish between connected components of the moduli space of smooth projective surfaces, see also \cite{FrTe}.

Moreover, the Zariski-Lefschetz hyperplane section theorem (see \cite{milnor})
states that
$$\pi_1 (\C\P^N \setminus S) \cong \pi_1 (H \setminus (H \cap S)),$$
where $S$ is an hypersurface and $H$ is a generic 2-plane.
Since $H \cap S$ is a plane curve, the fundamental groups of complements of curves
can be used also for computing the fundamental groups of complements of hypersurfaces in $\C\P^N$.

A different need for fundamental groups' computations arises in the search for more examples of Zariski pairs \cite{Z1,Z2}. A pair of plane
curves is called {\it a Zariski pair} if they have the
same combinatorics (to be exact: there is a degree-preserving bijection between the set of irreducible
components of the two curves $C_1,C_2$, and there exist regular neighbourhoods of the curves $T(C_1),T(C_2)$
such that the pairs $(T(C_1),C_1),(T(C_2),C_2)$ are homeomorphic and the homeomorphism respects the bijection
above \cite{AB-CR}), but their complements in $\P^2$ are not homeomorphic.
For a survey, see \cite{ABCT}.

It is also interesting to explore new finite non-abelian groups which
serve as fundamental groups of complements of plane curves in general,
see for example \cite{Z1,AB,AB1,Deg}.

\medskip

An arrangement of lines in $\C^2$ is a union of copies of $\C^1$ in $\C^2$. Such an arrangement is called {\em real} if the defining equations of the lines can be written with real coefficients, and {\em complex} otherwise. Note that the intersection of
the affine part of a real arrangement with the natural copy of $\R^2$ in $\C^2$ is an
arrangement of lines in the real plane.

\begin{figure}[!ht]
\epsfysize 4cm
\epsfbox{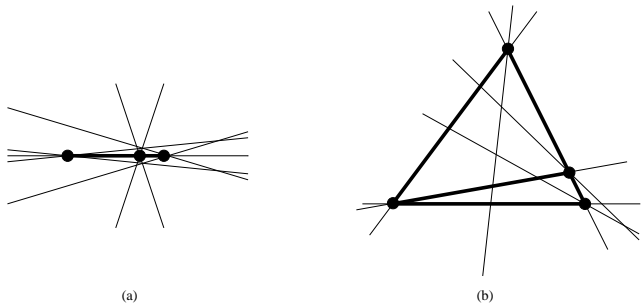}
\caption{Examples for $G(\mathcal L)$}\label{graph_GL}
\end{figure}

For real and complex line arrangements $\mathcal L$, Fan \cite{Fa2} defined a graph $G(\mathcal L)$ which is associated to its multiple points (i.e. points where more than two lines are intersected): Given a line arrangement $\mathcal L$, the graph $G(\mathcal L)$ of multiple points lies on $\mathcal L$. It consists of the multiple points of $\mathcal L$, with the
segments between the multiple points on lines which have at least two multiple points. Note that if the arrangement consists of three
multiple points on the same line, then $G(\mathcal L)$ has three vertices on the same line (see Figure \ref{graph_GL}(a)).
If two such lines happen to intersect in a simple point (i.e. a point where exactly two lines are intersected), it is ignored
(and the lines are not considered to meet in the graph theoretic sense). See another example in Figure \ref{graph_GL}(b) (note that this definition gives a graph different from the graph defined in \cite{JY}).

Fan \cite{Fa1,Fa2} proved some results concerning the projective fundamental group:

\begin{prop}[Fan]\label{Fan}
Let $\mathcal L$ be a complex arrangement of $n$ lines and\break $S=\{a_1, \dots, a_p\} $
be the set of all multiple points of $\mathcal L$.  Suppose that
$\beta (\mathcal L)=0$, where $\beta (\mathcal L)$ is the first Betti number of the graph $G(\mathcal L)$ (hence $\beta (\mathcal L)=0$
means that the graph $G(\mathcal L)$ has no cycles). Then:
$$\pi _1 (\C \P ^2 - \mathcal L) \cong \Z ^r \oplus \F _{m(a_1)-1} \oplus \cdots \oplus \F _{m(a_p)-1}$$
where $m(a_i)$ is the multiplicity of the intersection point $a_i$ and $r=n+p-1-m(a_1)- \cdots - m(a_p)$.
\end{prop}

In \cite{Ga,GaTe}, similar results were achieved for the affine and  projective
fundamental groups by different methods.

\medskip

Fan \cite{Fa2} has conjectured that the inverse implication is also correct, i.e. if the fundamental group $\pi _1 (\C \P ^2 - \mathcal L)$ can be written as a direct sum of free groups and infinite cyclic groups, then the graph $G(\mathcal L)$ has no cycles.

In an unpublished note, Fan \cite{Fa3} shows that if the fundamental group of the affine complement is a free group, then the arrangement consists of
parallel lines.

Recently, Eliyahu, Liberman, Schaps and Teicher \cite{ELST} proved Fan's conjecture completely.

\medskip

These results motivate the following definition:
\begin{defn}
Let $G$ be a fundamental group of the affine or projective complements of some line arrangement with $n$ lines. We say that $G$ has {\em a conjugation-free geometric presentation} if $G$ has a presentation with the following properties:
\begin{itemize}
\item In the affine case, the generators $\{ x_1,\dots, x_n \}$ are the meridians of lines at some far side of the arrangement, and therefore the number of generators is equal to $n$.
\item In the projective case, the generators are the meridians of lines at some far side of the arrangement except for one, and therefore the number of generators is equal to  $n - 1$.
\item In both cases, the relations are of the following type:
$$x_{i_k} x_{i_{k-1}} \cdots x_{i_1} = x_{i_{k-1}} \cdots x_{i_1} x_{i_k} = \cdots = x_{i_1} x_{i_k} \cdots x_{i_2},$$
where $\{ i_1,i_2, \dots , i_k \} \subseteq \{1, \dots, m \}$ is an increasing subsequence of indices,
where $m=n$ in the affine case and $m=n-1$ in the projective case. Note that for $k=2$ we get the usual commutator.
\end{itemize}

\end{defn}

Note that in usual geometric presentations of the fundamental group, most of the relations have conjugations (see Section \ref{MT}).

Based on the last definition, Fan's result yields that if the graph associated to the arrangement is acyclic, then the corresponding fundamental group has a conjugation-free geometric presentation.

\medskip

The following natural problem arises:
\begin{prob}
Which line arrangements have a fundamental group which has a  conjugation-free geometric presentation?
\end{prob}

The aim of this paper is to attack this problem.

\medskip

The importance of this family of arrangements is that the fundamental group can be read directly from the arrangement or equivalently from its incidence lattice (where the {\em incidence lattice} of an arrangement is the partially-ordered set of non-empty intersections of the lines, ordered by
inclusion, see \cite{OT}) without any computation. Hence, for this family of arrangements, the incidence lattice determines the fundamental group of the complement.

\medskip

%
%
%
%
%
%

\begin{figure}[!ht]
\epsfysize 4cm
\epsfbox{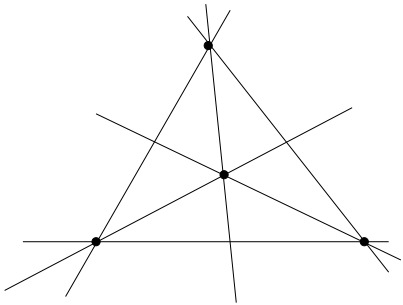}
\caption{Ceva arrangement}\label{ceva}
\end{figure}

We start with the easy fact that there exist arrangements whose fundamental groups have no conjugation-free geometric presentation: The fundamental group of the affine Ceva arrangement (also known as the {\it braid arrangement}, appears in Figure \ref{ceva}) has no conjugation-free geometric presentation. This fact was checked computationally by a package called {\it TESTISOM} \cite{HR}, which looks for isomorphisms (or proves a non-isomorphism) between two given finitely-presented group. Note that the Ceva arrangement is the minimal arrangement (with respect to the number of lines) with this property.

\medskip

Our main result is:

\begin{prop}
The fundamental group of following family of real arrangements have a conjugation-free geometric presentation: an arrangement $\mathcal L$, where $G(\mathcal L)$ is a union of disjoint cycles of any length, has no line with more than two multiple points, and the multiplicities of the multiple points are arbitrary.
\end{prop}

We also give the exact group structure (by means of a semi-direct product) of the fundamental group for an arrangement of $6$ lines whose graph is a cycle of length $3$ (i.e. a triangle), where all the multiple points are of multiplicity $3$:
\begin{prop}
Let $\mathcal L$ be the real arrangement of $6$ lines, whose graph consists of a cycle of length $3$, where all the multiple points are of multiplicity $3$. Moreover, it has no line with more than two multiple points. Then:
$$\pi_1(\C^2 -\mathcal L) \cong (\Z^2*\Z) \rtimes_{\alpha_3} \F_2 \rtimes_{\alpha_2} \F$$
where $*$ is the free product.
\end{prop}

\medskip

As mentioned above, for the family of arrangements with a conjugation-free geometric presentation of the fundamental group, the incidence lattice of the arrangement determines its fundamental group. There are some  well-known families of arrangements whose lattice determines the fundamental group of its complement - the families of {\it nice} arrangements (Jiang-Yau \cite{JY}) and {\it simple} arrangements (Wang-Yau \cite{WY}). It is interesting to study the relation between these families and the family of arrangements whose fundamental groups have conjugation-free geometric presentations, since for the latter family, the lattice determines the fundamental group of the complement too. We have the following remark:

\begin{rem}
The fundamental group of the arrangement $A_5$ (appears in Figure \ref{A_n}) has a conjugation-free geometric presentation (this fact was checked computationally), but this arrangement is neither nice nor simple.

\begin{figure}[!ht]
\epsfysize 4cm
\epsfbox{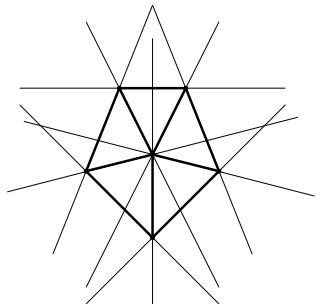}
\caption{The arrangement $A_5$}\label{A_n}
\end{figure}

\end{rem}

It will be interesting to find out whether our family of arrangements is broader than the family of simple arrangements, or whether there exists a simple arrangement whose fundamental group has no conjugation-free geometric presentation.

\medskip

\begin{rem}\label{rem_dehornoy}
It is worth to mention that conjugation-free geometric presentations are complemented positive presentations (defined by Dehornoy \cite{Deh1}, see also \cite{Deh2,Deh3}). Some initial computations show that in general conjugation-free geometric presentations are not complete (since the cube condition is not satisfied for some triples of generators). Nevertheless, we do think that there exist conjugation-free geometric presentations which are complete and hence have all the good properties induced by the completeness (see the survey \cite{Deh3}). We will discuss this subject in a different paper.
\end{rem}

\medskip

The paper is organized as follows. In Section \ref{MT}, we give a quick survey of the techniques we are using throughout the paper. In Section \ref{length_three}, we show that the fundamental group of a real arrangement whose graph has a unique cycle of length 3 has a conjugation-free geometric presentation. In this section, we also deal with the exact structure of the fundamental group of a real arrangement whose graph
consists of a cycle of length $3$, where all the multiple points have multiplicity $3$.
Section \ref{length_n} deals with the corresponding result for a real arrangement whose graph has a unique cycle of length $n$. We also generalize this result for the case of arrangements whose graphs are a union of disjoint cycles.


\section{The computation of the fundamental group}\label{MT}

In this section, we present the computation of the fundamental
group of the complement of real line arrangements. This
is based on the Moishezon-Teicher method
\cite{BGT1} and the van Kampen theorem \cite{vK}.
Some more presentations and algorithms can be found in \cite{arvola,CS,Ra,Sal}.

If the reader is familiar with this algorithm, he can skip this section.

\subsection{Wiring diagrams and Lefschetz pairs}

\begin{figure}[!ht]
\epsfysize 6cm
\epsfbox{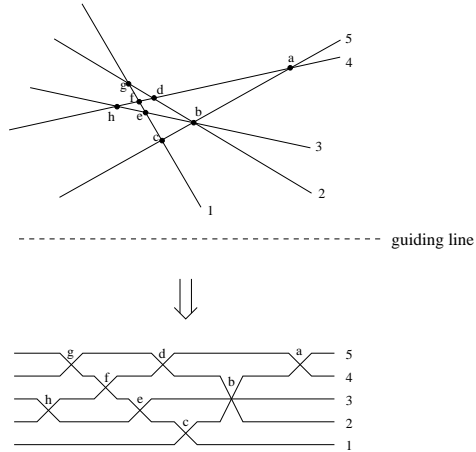}
\caption{From a line arrangement to a wiring diagram}\label{latowd}
\end{figure}

To an arrangement of $\ell$ lines in $\R^2$ one can associate a
\defin{wiring diagram} \cite{Go}, which holds the combinatorial
data of the arrangement and the position of the intersection
points. A wiring diagram is a collection of $\ell$ wires (where a
\defin{wire} in $\R ^2$ is a union of segments and rays,
homeomorphic to $\R$). The induced wiring diagram is constructed
by choosing a new line (called the \defin{guiding line}), which
avoids all the intersection points of the arrangement, such that
the projections of intersection points do not overlap. Then, the
$\ell$ wires are generated as follows. Start at the '$\infty$'
end of the line with $\ell$ parallel rays, and for every
projection of an intersection point, make the corresponding
switch in the rays, as in Figure \ref{latowd}.

To a wiring diagram, one can associate a list of \defin{Lefschetz pairs}.
Any pair of this list corresponds to one of the
intersection points, and holds the smallest
and the largest indices of the wires intersected at this point,
numerated locally near the intersection point (see \cite{BGT1} and \cite{GaTe}).

\begin{figure}[!ht]
\epsfysize 3cm
\epsfbox{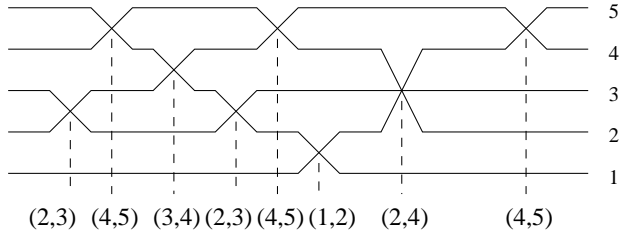}
\caption{Computing Lefschetz pairs for a wiring diagram}\label{wdtolp}
\end{figure}

For example, in the wiring diagram of Figure \ref{wdtolp},
the list of Lefschetz pairs is (we pass on the intersection points from right to left):
$$( \Lpair{4}{5},\Lpair{2}{4},\Lpair{1}{2},\Lpair{4}{5},\Lpair{2}{3},\Lpair{3}{4},\Lpair{4}{5},\Lpair{2}{3} ).$$

\subsection{Braid monodromy computation}
Let $D$ be a closed disk in $\R^2$, $K \sbs \Int(D)$ a set of $\ell$ points, and
$u \in\partial D$.
Let $\cB$ be the group of all diffeomorphisms
$\be: D\ra D$ such that $\be|_{\partial D}$ is the identity and $\be(K)=K$.
The action of such $\be$ on the disk applies to paths in $D$, which induces
an automorphism on
$\pi_1(D-K,u)$.
The \defin{braid group},
$B_{\ell}[D,K]$, is the group $\cB$ modulo the subgroup of
diffeomorphisms inducing the trivial automorphism on
$\pi_1(D-K,u)$. An element of
$B_{\ell} [D,K]$ is called a \defin{braid}.
For simplicity,
we will assume that $D = \set{ z \in \C \suchthat \power{z -
\frac{\ell+1}{2}} \leq \frac{\ell+1}{2} }$, and that $K = \set{ 1, 2, \dots , \ell } \sbs
D$.

\begin{figure}[!ht]
\epsfysize 4cm
\epsfbox{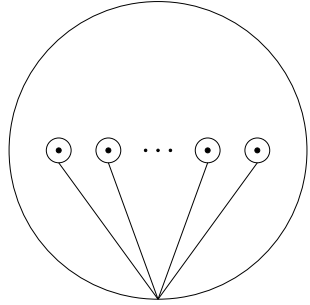}
\caption{A g-base}\label{fig_1}
\end{figure}

Choose a point $u_0 \in D$ (for convenience we choose it to be
below the real line). The group $\pi_1(D-K,u_0)$ is freely
generated by $\Gamma_1,\dots,\Gamma_\ell$, where $\Gamma_i$ is a
loop starting and ending at $u_0$, enveloping the \th{i} point in
$K$. The set $\set{\Gamma_1,\dots,\Gamma_\ell}$ is called a
\defin{geometric base} or \defin{g-base} of $\pi_1(D-K,u_0)$ (see Figure \ref{fig_1}).

Let $a=(\Lpair{a_1}{b_1},\dots,\Lpair{a_p}{b_p})$ be a list of
\LPs\ associated to a real line arrangement ${\mathcal L}$ with $\ell$ lines. The
\fg\ of the complement of the arrangement is a quotient group of
$\pi_1(D-K,u_0)$. There are $p$ sets of relations, one for every intersection point.
In each point, we will compute an object called a
\defin{skeleton}, from which the relation is computed.

\begin{figure}[!ht]
\epsfysize 0.7cm
\epsfbox{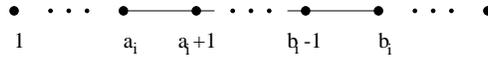}
\caption{The initial skeleton}\label{fig3}
\end{figure}

In order to compute the skeleton $s_i$ associated to the
\th{i} intersection point,
we start with an \defin{initial skeleton} corresponding to the  \th{i} Lefshetz pair $\Lpair{a_i}{b_i}$ which is presented in
Figure \ref{fig3}, in which the points correspond to the lines of
the arrangement and we connect by segments adjacent points which correspond to a local numeration of lines
passing through the intersection point.

To the initial skeleton, we apply the Lefschetz pairs
$\Lpair{a_{i-1}}{b_{i-1}},
\cdots,\Lpair{a_1}{b_1}$. A Lefschetz pair  $\Lpair{a_j}{b_j}$
acts by rotating the region from $a_j$ to $b_j$ by $180^{\circ}$
counterclockwise without affecting any other points.

\begin{figure}[!ht]
\epsfxsize 4.5cm
\epsfbox{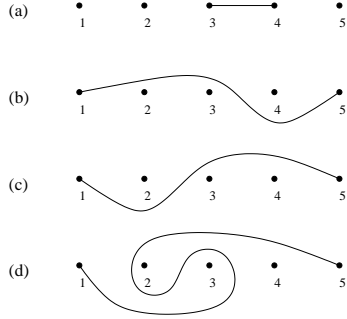}
\caption{An example for the computation of the braid monodromy}\label{fig4}
\end{figure}

For example, consider the list $\ba =
(\Lpair{2}{3},\Lpair{2}{4},\Lpair{4}{5},\Lpair{1}{3},\Lpair{3}{4})$.
Let us compute the skeleton associated to the 5th point.
The initial skeleton for $\Lpair{3}{4}$ is given in Figure
\ref{fig4}(a). By applying $\Lpair{1}{3}$ and then
$\Lpair{4}{5}$, we get the skeleton in Figure \ref{fig4}(b).
Then, applying  $\Lpair{2}{4}$ yields the skeleton in Figure
\ref{fig4}(c), and finally by
acting with  $\Lpair{2}{3}$ we get the final skeleton in Figure
\ref{fig4}(d).

\subsection{Inducing the presentation}

From the final skeletons we compute the relations, as follows.
We first explain the case when $\Lpair{c_i}{d_i}$ corresponds to a
simple point, \ie\ $d_i - c_i =1$. Then the skeleton is a path
connecting two points.

Let $D$ be a disk circumscribing the skeleton, and let $K$ be the set of points.
 Choose an arbitrary point on the path and 'pull' it down,
splitting the path into two parts, which are connected in one end
to $u_0 \in \partial D$ and in the other to the two end points of the path in $K$.

The loops associated to these two paths are elements in the group
$\pi_1 (D-K,u_0)$, and we call them $a_1$ and $a_2$. The corresponding
elements commute in the \fg\ of the arrangement's complement.

\begin{figure}[!ht]
\epsfysize 6cm
\epsfbox{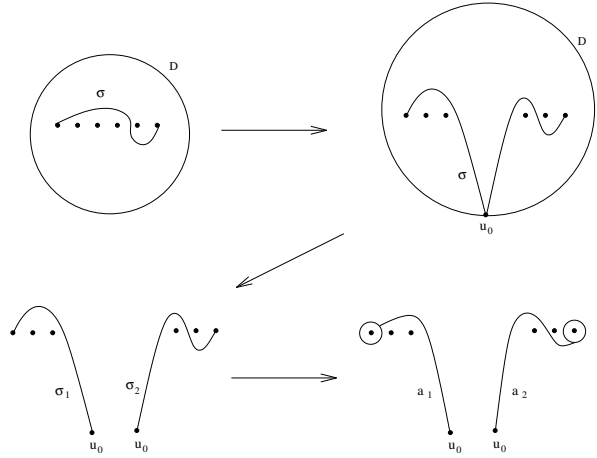}
\caption{Computation of $a_1,a_2$ for a simple intersection point}\label{av_bv}
\end{figure}

Figure \ref{av_bv} illustrates this procedure.

Now we show how to write $a_1$ and $a_2$ as words in the
generators
$\set{\Ga_1, \dots, \Ga_\ell}$ of $\pi _1(D-K,u_0)$. We start
with the generator corresponding to the end point of
$a_1$ (or $a_2$), and conjugate it as we move along $a_1$ (or
$a_2$) from its end point on $K$ to $u_0$ as follows: for every
point $i \in K$ which we pass from above, we conjugate by $\Ga_i$
when moving from left to right, and by $\Ga_i^{-1}$ when moving
from right to left.

For example, in Figure \ref{av_bv},
$$a_1 = \Ga_3 \Ga_2 \Ga_1 \Ga_2^{-1} \Ga_3^{-1}, \quad a_2 = \Ga_4 ^{-1} \Ga_6 \Ga_4,$$
and so the induced relation is:
$$\Ga_3 \Ga_2 \Ga_1 \Ga_2^{-1} \Ga_3^{-1} \cdot \Ga_4 ^{-1} \Ga_6 \Ga_4 =\Ga_4 ^{-1} \Ga_6 \Ga_4 \cdot \Ga_3 \Ga_2 \Ga_1 \Ga_2^{-1} \Ga_3^{-1}.$$

One can check that the relation is independent of the point in
which the path is split.

\bigskip

For a multiple intersection point of multiplicity $k$, we compute the elements in the
group
$\pi_1 (D-K,u_0)$ in a similar way, but the induced relations are
 of the following type:
$$a _k a_{k-1} \cdots  a _1 =
 a _1 a _k \cdots  a _3 a _2 = \cdots =
 a _{k-1} a _{k-2} \cdots a _1 a _k.$$
We choose an
arbitrary point on the path and pull it down to
$u_0$.
For each of the $k$ end points of the skeleton, we generate the
loop associated to the path from
$u_0$ to that point, and translate this path to a word in
$\G_1,\dots,\G_\ell$ by the procedure described above.

\begin{figure}[!ht]
\epsfysize 6cm
\epsfbox{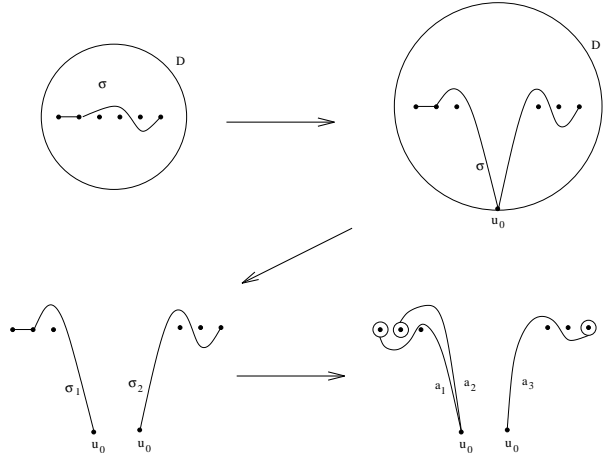}
\caption{Computation of $a_1,a_2,a_3$ for a multiple intersection point}\label{av_bv_mul}
\end{figure}

In the example given in Figure \ref{av_bv_mul}, we have:
$a_1 =
\G_3 \G_1 \G_3^{-1}$, $a_2 = \G_3 \G_2
\G_3^{-1}$ and $a_3 = \G_4^{-1} \G_6 \G_4$, so the relations are
\begin{eqnarray*}
\G_4^{-1} \G_6 \G_4 \cdot \G_3 \G_2 \G_3^{-1} \cdot \G_3 \G_1
\G_3^{-1} & = & \G_3 \G_1 \G_3^{-1} \cdot \G_4^{-1} \G_6 \G_4 \cdot
\G_3 \G_2 \G_3^{-1}\\
& = & \G_3 \G_2 \G_3^{-1} \cdot \G_3 \G_1
\G_3^{-1} \cdot \G_4^{-1} \G_6 \G_4.
\end{eqnarray*}


\section{An arrangement whose graph has a unique cycle of length $3$}\label{length_three}

In this section, we prove the following proposition:

\begin{prop}\label{triangle-prop}
The fundamental group of a real affine arrangement without parallel lines, whose graph which has a unique cycle of length $3$ and has no line with more than two multiple points, has a conjugation-free geometric presentation.
\end{prop}

In the first subsection we present the proof of Proposition \ref{triangle-prop}. The second subsection will be devoted to
studying the group structure of the fundamental group of the simplest arrangement of this family.

\subsection{Proof of Proposition \ref{triangle-prop}}
For simplicity, we will assume that all the multiple points have the same multiplicity $n+1$, but  the same argument will work even
if the multiplicities are not equal.

\medskip

\begin{figure}[!ht]
\epsfysize 9cm
\epsfbox{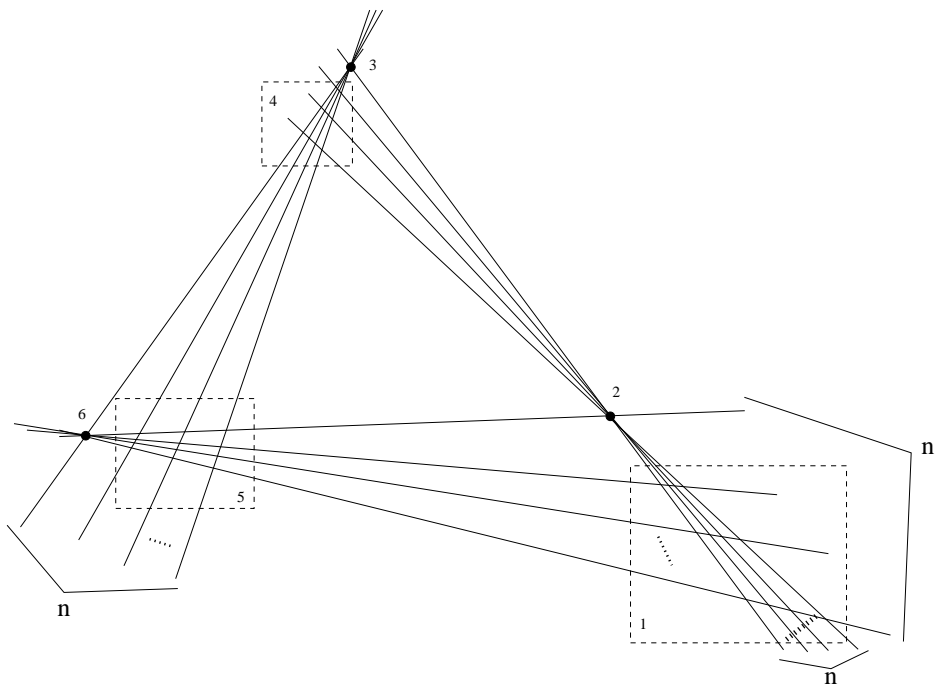}
\caption{The drawing of the arrangement with a cycle of length 3}\label{arrange_mult_3}
\end{figure}

\begin{figure}[!ht]
\epsfysize 3cm
\epsfbox{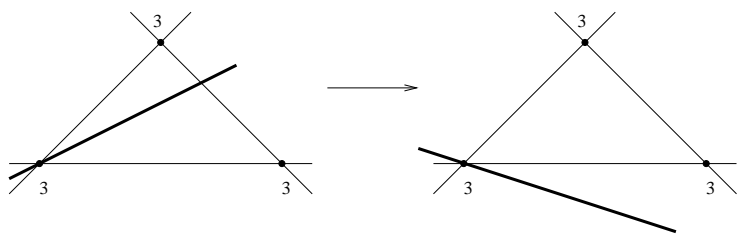}
\caption{Moving a line that participates in only one multiple point over a different line}\label{triangle-line}
\end{figure}

By rotations and translations, one can assume that we have a drawing of an  arrangement which has a unique cycle of multiple points of length 3 and has no line with more than two multiple points, as in Figure \ref{arrange_mult_3}. We can assume it due to the following reasons: First, one can rotate a line that participates in only one multiple point as long as it does not unite with a different line (by Results 4.8 and 4.13 of \cite{GTV}). Second, moving a line that participates in only one multiple point over a different line (see Figure \ref{triangle-line}) is permitted in the case of a triangle due to a result of Fan \cite{Fa2} that the family of configurations with $6$ lines and three triple points is connected by a finite sequence of smooth equisingular deformations.

Each of the blocks 1,4,5 contains simple intersection points of two pencils.

In block 1, one can assume that all the intersections of any horizontal line are adjacent, without intervening points from the third pencil.
In blocks 4 and 5, one can assume that all the intersections of any vertical line are adjacent
(in block 4, the vertical lines are those with positive slopes).
Moreover, all the intersection points of block 5 are to the left of all the intersection points of block 4.
Hence, we get the list of Lefschetz pairs as in Table \ref{tab1} (we put a double line to separate between the pairs related to different blocks).

\begin{tiny}
\begin{table}
$$\begin{array}{|c|c||c||c|c|}
\hline
j & {\rm Lefschetz\ pairs} &  &  j & {\rm Lefschetz\ pairs} \\
\hline\hline
1 & [ n,n+1 ] &  &  2n(n-1)+3 & [n,n+1] \\
2 & [ n-1,n ] &  &  2n(n-1)+4 & [n-1,n] \\
\vdots & \vdots  &  & \vdots & \vdots \\
n & [ 1,2 ] &  & n(2n-1)+2 & [1,2] \\
\hline
n+1 & [ n+1,n+2 ] & & n(2n-1)+3 & [n+1,n+2] \\
n+2 & [ n,n+1 ] & &  n(2n-1)+3 & [n,n+1] \\
\vdots & \vdots &  & \vdots & \vdots \\
2n & [ 2,3 ] &  & n(2n)+2 & [2,3] \\
\hline
\vdots & \vdots & & \vdots & \vdots\\
\hline
(n-2)n+1 & [ 2n-2,2n-1 ] &  & (3n-1)(n-1)+3 & [2n-2,2n-1]  \\
(n-2)n+2 & [ 2n-3,2n-2 ] &  & (3n-1)(n-1)+4 & [2n-3,2n-2] \\
\vdots & \vdots  &  & \vdots & \vdots\\
(n-1)n & [ n-1,n ] &   & 3n(n-1)+2 & [n-1,n] \\
\hline\hline
(n-1)n+1 & [ n,2n ] &  & 3n(n-1)+3 & [n,2n] \\
\hline\hline
(n-1)n+2 & [ 2n,3n ] & & & \\
\hline\hline
(n-1)n+3 & [ 2n-1,2n ] &  & &\\
(n-1)n+4 & [ 2n-2,2n-1 ] & & & \\
\vdots & \vdots  &  & &\\
(n-1)(n+1)+2 & [ n+1,n+2 ] & & & \\
\hline
(n-1)(n+1)+3 & [ 2n,2n+1 ] & & & \\
(n-1)(n+1)+4 & [ 2n-1,2n ] & & & \\
\vdots & \vdots &  & &\\
(n-1)(n+2)+2 & [ n+2,n+3 ] & & &  \\
\hline
\vdots & \vdots & & & \\
\hline
(2n-1)(n-1)+3 & [ 3n-2,3n-1 ] &  & &\\
(2n-1)(n-1)+4 & [ 3n-3,3n-2 ] &  & &\\
\vdots & \vdots  &  & &\\
2n(n-1)+2 & [ 2n,2n+1 ] &  & & \\
\hline
\end{array}
$$
\caption{List of Lefschetz pairs}\label{tab1}
\end{table}
\end{tiny}

\begin{figure}[!ht]
\epsfxsize 10cm
\epsfbox{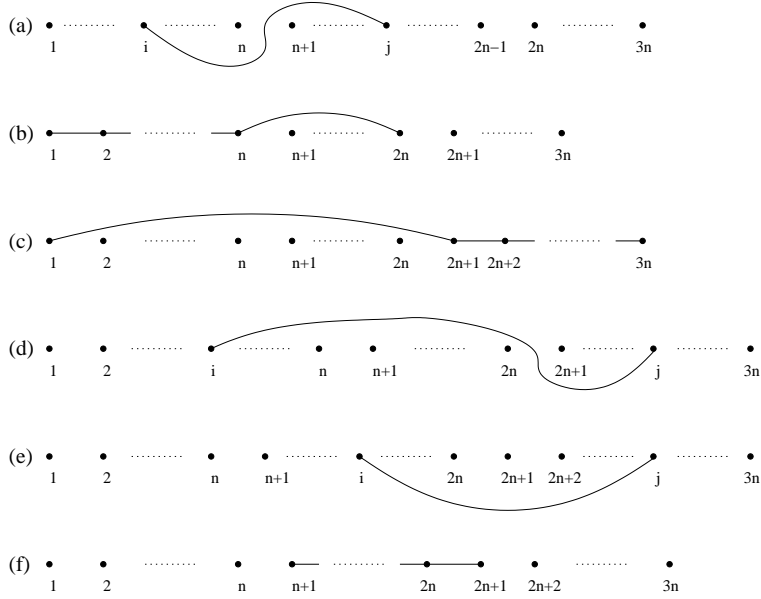}
\caption{The skeletons of the braid monodromy}\label{braid1-6}
\end{figure}

By the Moishezon-Teicher algorithm (see Section \ref{MT}), we get the following skeletons:
\begin{itemize}
\item For point $k$, where $1 \leq k \leq n(n-1)$, the corresponding final skeleton appears in
Figure \ref{braid1-6}(a), where $1 \leq i \leq n$ and $n+1 \leq j \leq 2n-1$.

\item For point $n(n-1)+1$, the corresponding final skeleton appears in
Figure \ref{braid1-6}(b).

\item For point $n(n-1)+2$, the corresponding final skeleton appears in
Figure \ref{braid1-6}(c).

\item For point $k$, where $n(n-1)+3 \leq k \leq 2n(n-1)+2$, the corresponding final skeleton
appears in Figure \ref{braid1-6}(d), where $2 \leq i \leq n$ and $2n+1 \leq j \leq 3n$.

\item For point $k$, where $2n(n-1)+3 \leq k \leq 3n(n-1)+2$, the corresponding final skeleton
appears in Figure \ref{braid1-6}(e), where $n+1 \leq i \leq 2n$ and $2n+2 \leq j \leq 3n$.

\item For point $3n(n-1)+3$, the corresponding final skeleton appears in
Figure \ref{braid1-6}(f).
\end{itemize}

Before we proceed to the presentation of the fundamental group, we introduce one notation: instead of writing the relations (where $a_i$ are words in a group):
$$a_n a_{n-1} \cdots a_1 = a_{n-1} \cdots a_1 a_n = \cdots =a_1 a_n \cdots a_2,$$
we will sometimes write: $[a_1,a_2 , \dots , a_n]$.

By the van Kampen theorem (see Section \ref{MT}), we get the following presentation of the fundamental group of the line arrangement's complement:

\medskip

\noindent
Generators: $\{x_1,x_2, \dots, x_{3n} \}$\\
Relations: \\
\begin{enumerate}
\item $[x_i, x_{n+1}^{-1} \cdots x_{j-1}^{-1} x_j x_{j-1} \cdots x_{n+1}]=e$, where $1 \leq i \leq n$ and\break $n+1 \leq j \leq 2n-1$.
\item $[x_1, x_2, \dots, x_n, x_{n+1}^{-1} \cdots x_{2n-1}^{-1} x_{2n} x_{2n-1} \cdots x_{n+1}]$.
\item $[x_{2n} x_{2n-1} \cdots x_2 x_1 x_2^{-1} \cdots x_{2n-1}^{-1} x_{2n}^{-1},x_{2n+1}, \dots, x_{3n}]$.
\item $[x_{2n} x_{2n-1} \cdots x_{i+1} x_i x_{i+1}^{-1} \cdots x_{2n-1}^{-1} x_{2n}^{-1}, x_j]=e$ where $2 \leq i \leq n$ and $2n+1 \leq j \leq 3n$.
\item $[x_i,x_j]=e$ where $n+1 \leq i \leq 2n$ and $2n+2 \leq j \leq 3n$.
\item $[x_{n+1}, \dots, x_{2n}, x_{2n+1}]$.
\end{enumerate}

\bigskip

Now, we show that all the conjugations in relations (1),(2),(3) and (4) can be simplified. We start with relations (1), and then relations (2).
We continue to relations (4) and we finish with relations (3).

We start with the first set of relations: for $j=n+1$, we get that for all $1 \leq i \leq n$ we have: $[x_i,x_{n+1}]=e$.
Now, we proceed to $j=n+2$. For $i=n$, we get: $[x_n,x_{n+1}^{-1} x_{n+2} x_{n+1}]=e$. By the relation $[x_n,x_{n+1}]=e$, it  is simplified to $[x_n,x_{n+2}]=e$. In this way, we get that for $j=n+2$, we have: $[x_i,x_{n+2}]=e$ for
$1 \leq i \leq n$.

By increasing $j$ one by one, we get that all the conjugations disappear and we get $[x_i, x_j]=e$, where $1 \leq i \leq n$ and $n+1 \leq j \leq 2n-1$, as needed.

\medskip

Relations (2) can be written as:
$$x_{n+1}^{-1} \cdots x_{2n-1}^{-1} x_{2n} x_{2n-1} \cdots x_{n+1} x_n \cdots x_1=$$
$$=x_1 x_{n+1}^{-1} \cdots x_{2n-1}^{-1} x_{2n} x_{2n-1} \cdots x_{n+1} x_n \cdots x_2 =$$
$$= \cdots = x_n \cdots x_1 x_{n+1}^{-1} \cdots x_{2n-1}^{-1} x_{2n} x_{2n-1} \cdots x_{n+1}$$
By the simplified version of relations (1), we can omit all the generators $x_{n+1}, \dots, x_{2n-1}$.
Hence we get:
$$ x_{2n} x_n \cdots x_1=x_1  x_{2n} x_n \cdots x_2 = \cdots = x_n \cdots x_1  x_{2n}$$
as needed.

\medskip

We proceed to relations (4). We start with $j=2n+1$. Taking $i=n$, we get:
$$[x_{2n} x_{2n-1} \cdots x_{n+1} x_n x_{n+1}^{-1} \cdots x_{2n-1}^{-1} x_{2n}^{-1}, x_{2n+1}]=e.$$
By relations (6), we have: $x_{2n} x_{2n-1} \cdots x_{n+1} x_n = x_n x_{2n} x_{2n-1} \cdots x_{n+1}$, and hence we get:
$$[x_n, x_{2n+1}]=e.$$
For $i=n-1$, we get:
$$[x_{2n} x_{2n-1} \cdots x_{n+1} x_n x_{n-1} x_n^{-1} x_{n+1}^{-1} \cdots x_{2n-1}^{-1} x_{2n}^{-1}, x_{2n+1}]=e.$$
By relations (6) again and the simplified version of relations (1), we get:
$$[x_n x_{n-1} x_n^{-1}, x_{2n+1}]=e.$$
Using the simplified relation $[x_n, x_{2n+1}]=e$, we get $[x_{n-1}, x_{2n+1}]=e$.

In the same way, we get that for $j=2n+1$ and $1\leq i \leq n$, we have: $[x_i, x_{2n+1}]=e$.

We continue to $j=2n+2$. Taking $i=n$, we have:
$$[x_{2n} x_{2n-1} \cdots x_{n+1} x_n x_{n+1}^{-1} \cdots x_{2n-1}^{-1} x_{2n}^{-1}, x_{2n+2}]=e.$$
By the simplified version of relations (1), we can omit all the generators $x_{n+1}, \dots , x_{2n-1}$. Hence we get:
$$[x_{2n} x_n x_{2n}^{-1}, x_{2n+2}]=e.$$
By relations (5), we can omit $x_{2n}$ too, and therefore: $[ x_n , x_{2n+2}]=e$.

For $i=n-1$, we have:
$$[x_{2n} x_{2n-1} \cdots x_{n+1} x_n x_{n-1} x_n^{-1} x_{n+1}^{-1} \cdots x_{2n-1}^{-1} x_{2n}^{-1}, x_{2n+2}]=e.$$
By the simplified version of relations (1), we can omit all the generators $x_{n+1}, \dots , x_{2n-1}$. Hence we get:
$$[x_{2n} x_n x_{n-1} x_n^{-1} x_{2n}^{-1}, x_{2n+2}]=e.$$
By relations (5), we can omit $x_{2n}$ too, and therefore:
$$[x_n x_{n-1} x_n^{-1} , x_{2n+2}]=e.$$
By $[ x_n , x_{2n+2}]=e$, we get $[ x_{n-1} , x_{2n+2}]=e$.
In the same way, we get that for $j=2n+2$ and $1\leq i \leq n$, we get: $[x_i, x_{2n+2}]=e$.

In the same way, by increasing $j$ one by one, we will get
that for all $2n+3 \leq j \leq 3n$ and $1\leq i \leq n$, we get: $[x_i, x_j]=e$ as needed.

\medskip

Relations (3) can be written:
$$x_{3n} \cdots x_{2n+1} x_{2n} x_{2n-1} \cdots x_2 x_1 x_2^{-1} \cdots x_{2n-1}^{-1} x_{2n}^{-1}=$$
$$=x_{3n-1} \cdots x_{2n+1} x_{2n} x_{2n-1} \cdots x_2 x_1 x_2^{-1} \cdots x_{2n-1}^{-1} x_{2n}^{-1} x_{3n}= \cdots =$$
$$=x_{2n} x_{2n-1} \cdots x_2 x_1 x_2^{-1} \cdots x_{2n-1}^{-1} x_{2n}^{-1} x_{3n} \cdots x_{2n+1}$$
By relations (1), we can omit the generators $x_{n+1}, \dots, x_{2n-1}$, so we get:
$$x_{3n} \cdots x_{2n+1} x_{2n} x_n \cdots x_2 x_1 x_2^{-1} \cdots x_n^{-1} x_{2n}^{-1}=$$
$$=x_{3n-1} \cdots x_{2n+1} x_{2n} x_{n} \cdots x_2 x_1 x_2^{-1} \cdots x_{n}^{-1} x_{2n}^{-1} x_{3n}= \cdots =$$
$$=x_{2n} x_{n} \cdots x_2 x_1 x_2^{-1} \cdots x_{n}^{-1} x_{2n}^{-1} x_{3n} \cdots x_{2n+1}$$
By relations (2), we can omit also the generators $x_{2},x_3, \dots, x_{n}, x_{2n}$ in order to get:
$$x_{3n} \cdots x_{2n+1} x_1 = x_{3n-1} \cdots x_{2n+1} x_1 x_{3n}= \cdots = x_1 x_{3n} \cdots x_{2n+1}.$$

Hence, we get the following simplified presentation:

\noindent
Generators: $\{x_1,x_2, \dots, x_{3n} \}$\\
Relations:

\medskip

\begin{enumerate}
\item $[x_i, x_j]=e$, where $1 \leq i \leq n$ and $n+1 \leq j \leq 2n-1$.
\item $[x_1, x_2, \dots, x_n, x_{2n}]$.
\item $[x_1,x_{2n+1}, \dots, x_{3n}]$.
\item $[x_i, x_j]=e$ where $2 \leq i \leq n$ and $2n+1 \leq j \leq 3n$.
\item $[x_i,x_j]=e$ where $n+1 \leq i \leq 2n$ and $2n+2 \leq j \leq 3n$.
\item $[x_{n+1}, \dots, x_{2n}, x_{2n+1}]$.
\end{enumerate}
Therefore, we have a conjugation-free geometric presentation, and hence we are done. \hfill \qed

\subsection{The structure of the fundamental group of the simplest case of this family}

Cohen and Suciu \cite{CS-Hom} give the following presentation of $\F_3 \rtimes_{\alpha_3} \F_2 \rtimes_{\alpha_2} \F_1$, which is known \cite{FN} to be the fundamental group of the complement of the affine Ceva arrangement (see Figure \ref{ceva}):
$$\F_1=\langle u \rangle, \qquad \F_2=\langle t,s \rangle, \qquad \F_3=\langle x,y,z \rangle$$
The actions of the automorphisms $\alpha_2$ and $\alpha_3$ are defined as follows:

\medskip

\noindent
$(\alpha_2(u))(t)=sts^{-1}, \qquad \ \ (\alpha_2(u))(s)=stst^{-1}s^{-1}, $\\
$(\alpha_2(u))(x)=x,\qquad \qquad (\alpha_2(u))(y)=zyz^{-1},\qquad (\alpha_2(u))(z)=zyzy^{-1}z^{-1}$

\medskip

\noindent
$(\alpha_3(s))(x)=zxz^{-1},\qquad (\alpha_3(s))(y)=zxz^{-1}x^{-1}yxzx^{-1}z^{-1},$ \\ $(\alpha_3(s))(z)=zxzx^{-1}z^{-1}$

\medskip

\noindent
$(\alpha_3(t))(x)=yxy^{-1},\qquad (\alpha_3(t))(y)=yxyx^{-1}y^{-1},$ \\ $(\alpha_3(t))(z)=z$

\medskip

Notice that if we rotate clockwise the lowest line in the affine Ceva arrangement (Figure \ref{ceva}), we get an arrangement $\mathcal L$ whose graph consists of a unique cycle of length $3$, where all the multiple points are of multiplicity $3$.

By a simple check, the effect of this rotation is the addition of the commutator relation $[x,z]=e$ to the presentation of the group. Hence, we get that the actions of the automorphisms $\alpha_2$ and $\alpha_3$  are changed as follows:

\medskip

\noindent
$(\alpha_2(u))(t)=sts^{-1}, \qquad \ \ (\alpha_2(u))(s)=stst^{-1}s^{-1}, $\\
$(\alpha_2(u))(x)=x,\qquad \qquad (\alpha_2(u))(y)=zyz^{-1},\qquad (\alpha_2(u))(z)=zyzy^{-1}z^{-1}$

\medskip

\noindent
$(\alpha_3(s))(x)=x,\qquad \qquad(\alpha_3(s))(y)=y, \qquad \qquad(\alpha_3(s))(z)=z$

\medskip

\noindent
$(\alpha_3(t))(x)=yxy^{-1},\qquad (\alpha_3(t))(y)=yxyx^{-1}y^{-1},$ \\ $(\alpha_3(t))(z)=z$

\medskip

This is the presentation of the group:  $(\Z^2*\Z) \rtimes_{\alpha_3} \F_2 \rtimes_{\alpha_2} \F$, where $*$ is the free product.

\medskip

To summarize, we get the following result:
\begin{prop}
Let $\mathcal L$ be the arrangement of $6$ lines without parallel lines whose graph is a unique cycle of length $3$, where all the multiple points are of multiplicity $3$. Then:
$$\pi_1(\C^2 -\mathcal L) \cong (\Z^2*\Z) \rtimes_{\alpha_3} \F_2 \rtimes_{\alpha_2} \F$$
\end{prop}

It is interesting to check how this proposition can be generalize to arrangements whose graphs are cycles of length $n>3$.


\section{An arrangement whose graph is a cycle of length $n$}\label{length_n}

In this section, we show that the fundamental group of a real affine arrangement whose graph is a unique cycle of any length and has no line with more than two multiple points, has a conjugation-free geometric presentation. At the end of this section, we generalize this result to arrangements whose graphs are unions of disjoint cycles.

\medskip

We start by investigating the case of a cycle of length $5$ and then we generalize it to any length.

\begin{figure}[!ht]
\centerline{\epsfxsize 12cm
\epsfbox{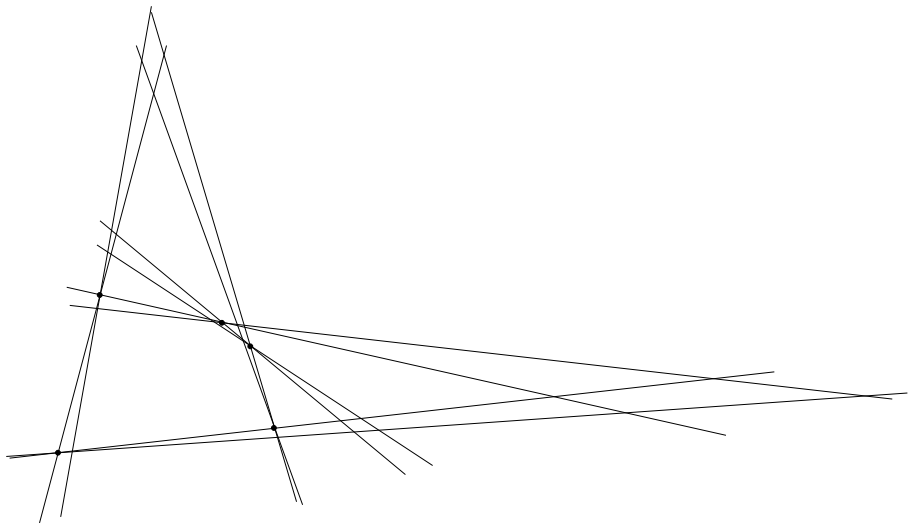}}
\caption{An arrangement whose graph is a cycle of length 5}\label{cycle_multiple5}
\end{figure}

In Figure \ref{cycle_multiple5}, we present a real arrangement whose graph is a cycle of $5$ multiple points (note that any real arrangement whose graph is a unique cycle of $5$ multiple points and has no line with more than two multiple points, can be transferred to this drawing by rotations, translations and equisingular deformations).

Based on Figure \ref{cycle_multiple5}, we get the list of Lefschetz  pairs presented in Table \ref{tab2}.

\begin{table}
$$\begin{array}{||c|c|c||c||c|c|c||c||c|c|c||}
\hline\hline
j & {\rm LP} & {\rm Mult} &\ \ & j & {\rm LP} & {\rm Mult} &\ \ & j & {\rm LP} & {\rm Mult} \\
\hline\hline
1 & [ 6,7 ]  & 2 & & 13 & [ 6,7 ] & 2 & & 25 & [ 8,9 ] & 2 \\
2 & [ 5,6 ]  & 2 & & 14 & [ 7,8 ] & 2 & & 26 & [ 6,7 ] & 2 \\
3 & [ 7,8 ]  & 2 & & 15 & [ 3,4 ] & 2 & & 27 & [ 5,6 ] & 2 \\
4 & [ 6,7 ]  & 2 & & 16 & [ 4,5 ] & 2 & & 28 & [ 7,8 ] & 2 \\
5 & [ 4,5 ]  & 2 & & 17 & [ 5,6 ] & 2 & & 29 & [ 6,7 ] & 2 \\
6 & [ 3,4 ]  & 2 & & 18 & [ 6,7 ] & 2 & & 30 & [ 4,6 ] & 3 \\
7 & [ 5,6 ]  & 2 & & 19 & [ 4,6 ] & 3 & & 31 & [ 3,4 ] & 2 \\
8 & [ 4,5 ]  & 2 & & 20 & [ 3,4 ] & 2 & & 32 & [ 4,5 ] & 2 \\
9 & [ 2,3 ]  & 2 & & 21 & [ 4,5 ] & 2 & & 33 & [ 2,3 ] & 2 \\
10 & [ 1,2 ] & 2 & & 22 & [ 8,9 ] & 2 & & 34 & [ 1,2 ] & 2 \\
11 & [ 2,4 ] & 3 & & 23 & [ 7,8 ] & 2 & & 35 & [ 2,4 ] & 3 \\
12 & [ 4,6 ] & 3 & & 24 & [ 9,10 ] & 2 & &  &  & \\
\hline\hline
\end{array}
$$
\caption{List of Lefschetz pairs of the arrangement in Figure \ref{cycle_multiple5}}\label{tab2}
\end{table}

\begin{figure}[!ht]
\epsfxsize 12cm
\epsfbox{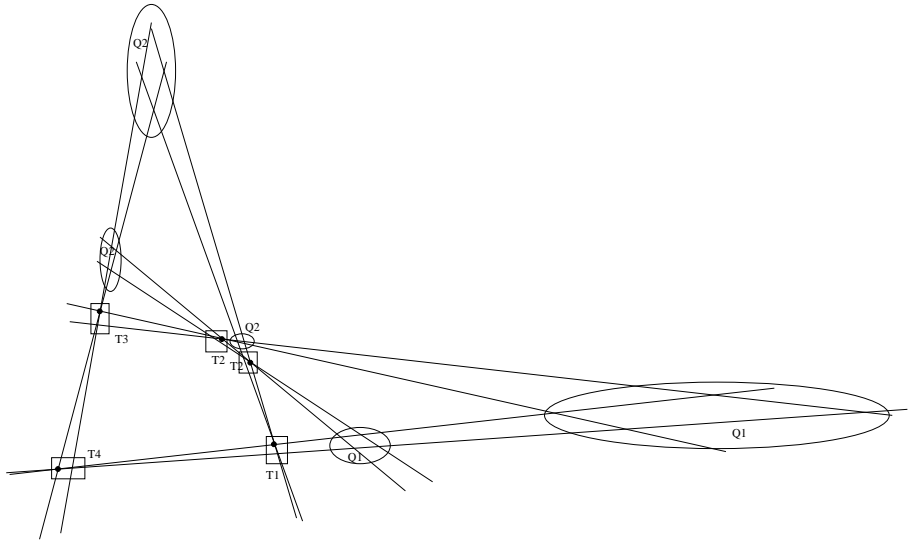}
\caption{The arrangement where the points are grouped into blocks}\label{cycle_multiple5_block}
\end{figure}

By the Moishezon-Teicher algorithm (see Section \ref{MT}), one can compute the skeletons of the braid monodromy. 
After the computation, one should notice that actually we can group the intersection points into blocks according to their braid monodromies (see Figure \ref{cycle_multiple5_block}), since the structure of the skeletons is similar.


Following this observation, we can deal with each block separately. So, we get the
following sets of skeletons:

\begin{itemize}

\item Quadruples of type Q1: see Figure \ref{quadruple_case1}(a) for $2 \leq i \leq 3$.

\begin{figure}[!ht]
\epsfxsize 11cm
\epsfbox{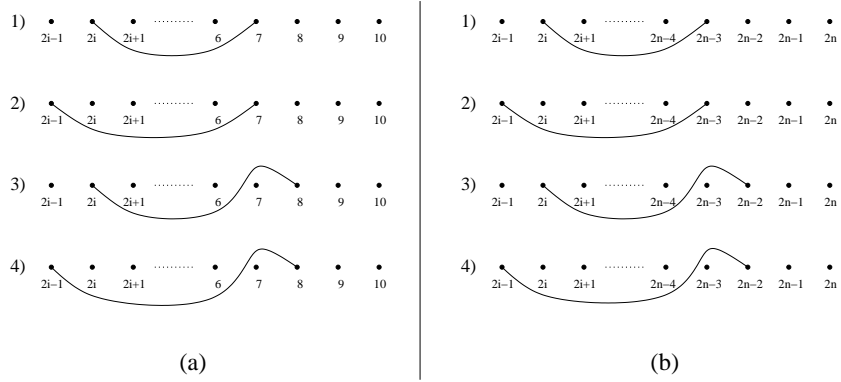}
\caption{Skeletons of quadruple Q1}\label{quadruple_case1}
\end{figure}

\item Quadruples of type Q2: see Figure \ref{quadruple_case2}(a) for $i,j \neq 4$, $|i-j|>1$, $(i,j)\neq(3,5)$.

\begin{figure}[!ht]
\epsfxsize 11cm
\epsfbox{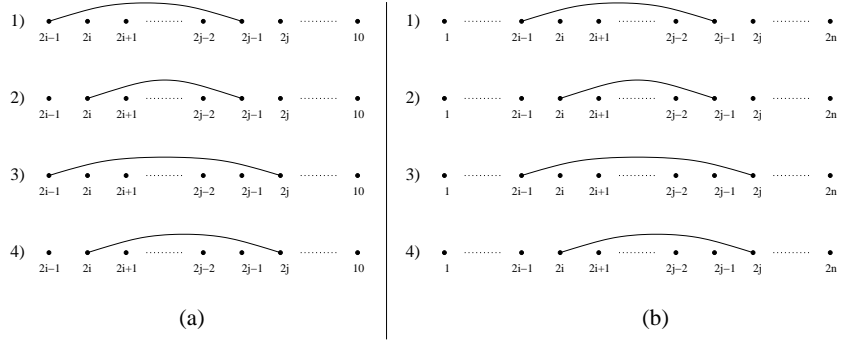}
\caption{Skeletons of quadruple Q2}\label{quadruple_case2}
\end{figure}

\item A triple of type T1: see Figure \ref{triple_case1}(a).

\begin{figure}[!ht]
\epsfxsize 11cm
\epsfbox{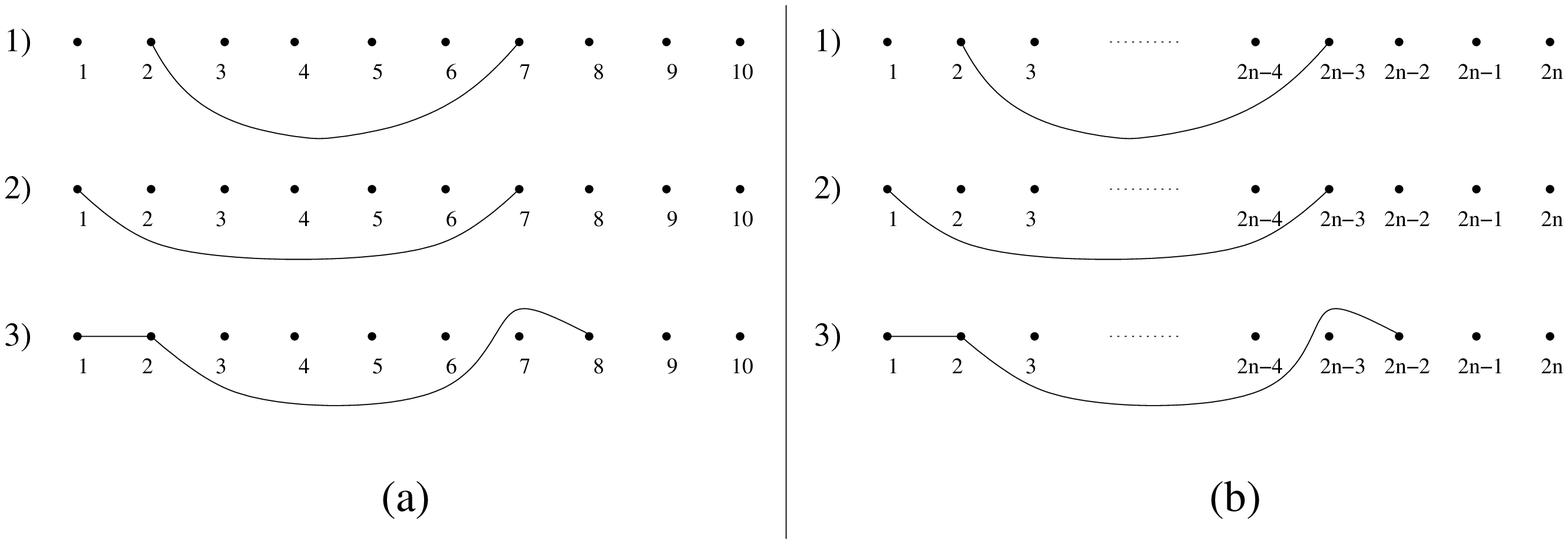}
\caption{Skeletons of triple T1}\label{triple_case1}
\end{figure}

\item Triples of type T2: see Figure \ref{triple_case2} for $1 \leq i \leq 2$.

\begin{figure}[!ht]
\epsfysize 4cm
\epsfbox{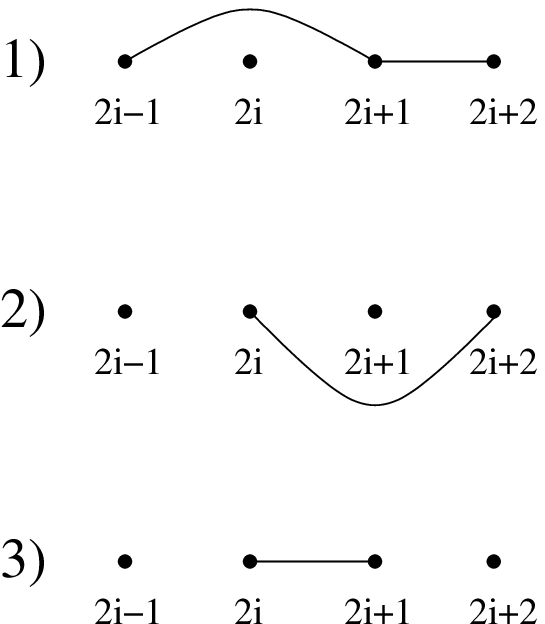}
\caption{Skeletons of triple T2}\label{triple_case2}
\end{figure}

\item A triple of type T3: see Figure \ref{triple_case3}(a).

\begin{figure}[!ht]
\epsfxsize 11cm
\epsfbox{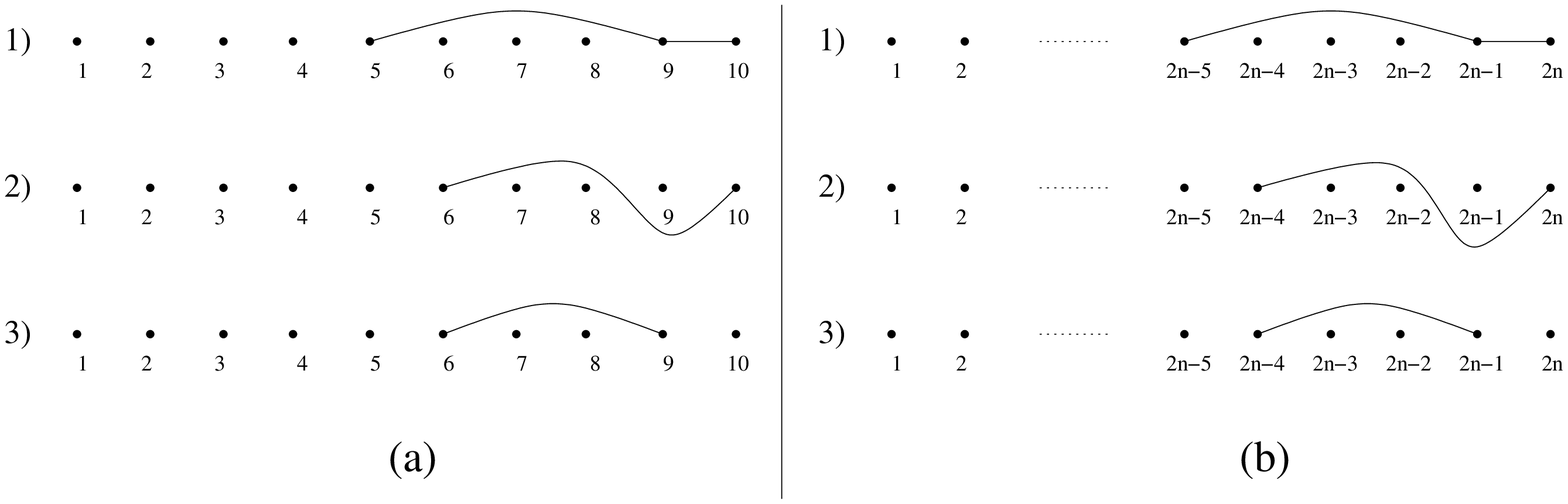}
\caption{Skeletons of triple T3}\label{triple_case3}
\end{figure}

\item A triple of type T4: see Figure \ref{triple_case4}(a).

\begin{figure}[!ht]
\epsfxsize 11cm
\epsfbox{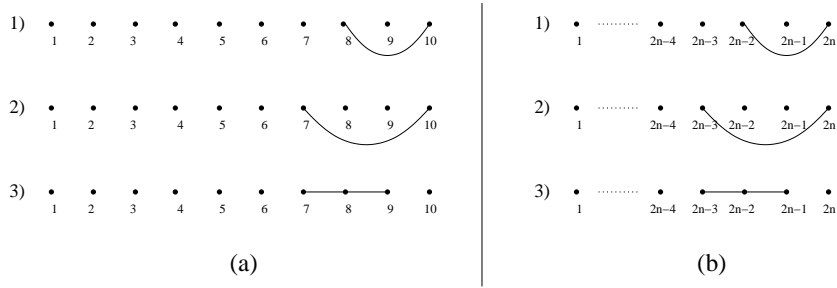}
\caption{Skeletons of triple T4}\label{triple_case4}
\end{figure}

\end{itemize}

\medskip

Now we pass to the general case. One can draw an arrangement of $2n$ lines whose graph is a unique cycle of length $n$ and has no line with more than two multiple points in a similar way to the way we have drawn the arrangement of $10$ lines whose graph is a cycle of length $5$. Hence, one can compute the braid monodromy of the general arrangement in blocks similar to what we have done in the case of $n=5$:

\medskip

\begin{itemize}

\item Quadruples of type Q1: for $2 \leq i \leq n-2$, see Figure \ref{quadruple_case1}(b).
\item Quadruples of type Q2: for $i,j \neq n-1$, $|i-j|>1$, $(i,j) \neq (n-2,n)$, see Figure \ref{quadruple_case2}(b).
\item A triple of type T1: see Figure \ref{triple_case1}(b).
\item Triples of type T2: for $1 \leq i \leq n-3$, see Figure \ref{triple_case2}.
\item A triple of type T3: see Figure \ref{triple_case3}(b).
\item A triple of type T4: see Figure \ref{triple_case4}(b).
\end{itemize}

By the van-Kampen theorem (see Section \ref{MT}), we get the following presentation of the fundamental group of the complement of the arrangement:

\medskip

\noindent
Generators: $\{ x_1, \dots, x_{2n} \}$\\
Relations:
\begin{itemize}

\item From quadruples of type Q1:

\begin{enumerate}
\item $[x_{2i},x_{2n-3}]=e$ where $2 \leq i \leq n-2$
\item $[x_{2i-1},x_{2n-3}]=e$ where $2 \leq i \leq n-2$
\item $[x_{2i},x_{2n-3}^{-1} x_{2n-2} x_{2n-3}]=e$ where $2 \leq i \leq n-2$
\item $[x_{2i-1},x_{2n-3}^{-1} x_{2n-2} x_{2n-3}]=e$ where $2 \leq i \leq n-2$
\end{enumerate}

\medskip

\item From quadruples of type Q2: for $i,j \neq n-1$, $|i-j|>1$, $(i,j) \neq (n-2,n)$:

\begin{enumerate}
\item $[x_{2i}, x_{2i+1}^{-1} \cdots x_{2j-1}^{-1} x_{2j} x_{2j-1} \cdots x_{2i+1}]=e$
\item $[x_{2i-1}, x_{2i}^{-1} x_{2i+1}^{-1} \cdots x_{2j-1}^{-1} x_{2j} x_{2j-1} \cdots x_{2i+1} x_{2i}]=e$
\item $[x_{2i}, x_{2i+1}^{-1} \cdots x_{2j-2}^{-1} x_{2j-1} x_{2j-2} \cdots x_{2i+1}]=e$
\item $[x_{2i-1}, x_{2i}^{-1} x_{2i+1}^{-1} \cdots x_{2j-2}^{-1} x_{2j-1} x_{2j-2} \cdots x_{2i+1} x_{2i}]=e$
\end{enumerate}

\medskip

\item From the triple of type T1:

\begin{enumerate}
\item $[x_2,x_{2n-3}]=e$
\item $[x_1,x_{2n-3}]=e$
\item $x_{2n-3}^{-1} x_{2n-2} x_{2n-3} x_2 x_1 = x_2 x_1 x_{2n-3}^{-1} x_{2n-2} x_{2n-3} = x_1 x_{2n-3}^{-1} x_{2n-2} x_{2n-3} x_2$
\end{enumerate}

\medskip

\item From triples of type T2:

\begin{enumerate}
\item $x_{2i+2} x_{2i+1} x_{2i} x_{2i-1} x_{2i}^{-1} = x_{2i+1} x_{2i} x_{2i-1} x_{2i}^{-1} x_{2i+2} = x_{2i} x_{2i-1} x_{2i}^{-1} x_{2i+2} x_{2i+1}$ where $1 \leq i \leq n-3$
\item $[x_{2i},x_{2i+2}]=e$ where $1 \leq i \leq n-3$
\item $[x_{2i},x_{2i+1}]=e$ where $1 \leq i \leq n-3$
\end{enumerate}

\medskip

\item From the triple of type T3:

\begin{enumerate}
\item $x_{2n} x_{2n-1} x_{2n-2} x_{2n-3} x_{2n-4} x_{2n-5} x_{2n-4}^{-1} x_{2n-3}^{-1} x_{2n-2}^{-1} =\\ = x_{2n-1} x_{2n-2} x_{2n-3} x_{2n-4} x_{2n-5} x_{2n-4}^{-1} x_{2n-3}^{-1} x_{2n-2}^{-1} x_{2n} =\\= x_{2n-2} x_{2n-3} x_{2n-4} x_{2n-5} x_{2n-4}^{-1} x_{2n-3}^{-1} x_{2n-2}^{-1} x_{2n} x_{2n-1}$
\item $[x_{2n-4}, x_{2n-3}^{-1} x_{2n-2}^{-1} x_{2n} x_{2n-2} x_{2n-3}]=e$
\item $[x_{2n-4}, x_{2n-3}^{-1} x_{2n-2}^{-1} x_{2n-1} x_{2n-2} x_{2n-3}]=e$
\end{enumerate}

\medskip

\item From the triple of type T4:

\begin{enumerate}
\item $[x_{2n-2}, x_{2n}]=e$
\item $[x_{2n-3}, x_{2n}]=e$
\item $x_{2n-1} x_{2n-2} x_{2n-3} = x_{2n-2} x_{2n-3} x_{2n-1} = x_{2n-3} x_{2n-1} x_{2n-2}$
\end{enumerate}

\end{itemize}

\medskip

We now show that all the conjugations can be simplified, and hence we have a  conjugation-free geometric presentation for the fundamental group.

We have conjugations in the relations coming from triples of points and quadruples of points.
We start with the relations which correspond to triples of points.

The conjugation in relation (3) of the triple of type T1 can be simplified using relations (1) and (2) of the triple of type T1.

The conjugation in relation (1) of triples of type T2 can be simplified using relations (2) and (3) of the corresponding triples of type T2.

The conjugation in relation (3) of the triple of type T3 can be simplified using relation (3) of the triple of type T4. The conjugation in relation (2) of the triple of type T3 can be simplified using relations (1) and (2) of the triple of type T4. The conjugation in relation (1) of the triple of type T3 can be simplified using relations (2)--(3) of the triple of type T3 and relations (1)--(3) of the triple of type T4.

\medskip

We continue to the relations induced by to the quadruples of type Q1.
By the first two relations of the quadruples of type Q1, one can easily simplify
the conjugations which appear in the last two relations of the quadruples of type Q1.
So we get that the relations correspond to the quadruples of type Q1 can be written without conjugations.

\medskip

Now, we pass to the relations correspond to the quadruples of type Q2.
We start with $i=1, j=3$: We have the following relations:
\begin{enumerate}
\item[(a)] $[x_2, x_3^{-1} x_4^{-1} x_5^{-1} x_6 x_5 x_4 x_3]=e$
\item[(b)] $[x_1, x_3^{-1} x_4^{-1} x_5^{-1} x_6 x_5 x_4 x_3]=e$
\item[(c)] $[x_1, x_3^{-1} x_4^{-1} x_5 x_4 x_3]=e$
\item[(d)] $[x_2, x_3^{-1} x_4^{-1} x_5 x_4 x_3]=e$
\end{enumerate}

By relations (2) and (3) of triple T2 (for $i=1$), we have the relations $[x_2,x_3]=e$ and $[x_2,x_4]=e$. Hence, the conjugations in relation (d) are canceled and we get $[x_2,x_5]=e$. By the same relations and the simplified version of relation (d), we get
the following relation from relation (a): $[x_2,x_6]=e$.

Substituting $i=1$ in relation (1) of triple T2 yields $$x_4 x_3 x_1 = x_3 x_1 x_4 = x_1 x_4 x_3.$$ By this relation, relation (c) becomes $[x_1,x_5]=e$, and relation (b) becomes $[x_1,x_6]=e$.

The same argument holds for any $i,j$, where $j-i=2$ and $j \leq n-2$. Hence, one can simplify the conjugations in these cases.

\medskip

Now, we pass to the case where $j-i=3$ and $j \leq n-2$. Let $i=1, j=4$. We have the following relations:
\begin{enumerate}
\item[(a')] $[x_2, x_3^{-1} x_4^{-1} x_5^{-1} x_6^{-1} x_7^{-1} x_8 x_7 x_6 x_5 x_4 x_3]=e$
\item[(b')] $[x_1, x_3^{-1} x_4^{-1} x_5^{-1} x_6^{-1} x_7^{-1} x_8 x_7 x_6 x_5 x_4 x_3]=e$
\item[(c')] $[x_1, x_3^{-1} x_4^{-1} x_5^{-1} x_6^{-1} x_7 x_6 x_5 x_4 x_3]=e$
\item[(d')] $[x_2, x_3^{-1} x_4^{-1} x_5^{-1} x_6^{-1} x_7 x_6 x_5 x_4 x_3]=e$
\end{enumerate}

By relations (2) and (3) of triple T2 (for $i=1$), we have the relations $[x_2,x_3]=e$ and $[x_2,x_4]=e$. Hence, relations (a') and (d') become:
\begin{enumerate}
\item[(a')] $[x_2, x_5^{-1} x_6^{-1} x_7^{-1} x_8 x_7 x_6 x_5]=e$
\item[(d')] $[x_2, x_5^{-1} x_6^{-1} x_7 x_6 x_5]=e$
\end{enumerate}
By relations (a) and (d) above, we get $[x_2,x_7]=e$, and therefore we also get $[x_2,x_8]=e$.

\medskip

Substituting $i=1$ in relation (1) of triple T2 yields $$x_4 x_3 x_1 = x_3 x_1 x_4 = x_1 x_4 x_3.$$ By this relation, relations (b') and (c')  become:
\begin{enumerate}
\item[(b')] $[x_1, x_5^{-1} x_6^{-1} x_7^{-1} x_8 x_7 x_6 x_5]=e$
\item[(c')] $[x_1, x_5^{-1} x_6^{-1} x_7 x_6 x_5]=e$
\end{enumerate}

By relations (b) and (c) above, we get  $[x_1, x_7]=e$ and hence $[x_1, x_8]=e$.

\medskip

It is easy to show by a simple induction that we can simplify the conjugations for any $i,j$, where $|j-i|>1$ and $j \leq n-2$.

\medskip

The remaining case is $j=n$. We start with $i=n-3$. We have the following relations:\\
\begin{tiny}
(a) $[x_{2n-6}, x_{2n-5}^{-1} x_{2n-4}^{-1} x_{2n-3}^{-1} x_{2n-2}^{-1} x_{2n-1}^{-1} x_{2n} x_{2n-1} x_{2n-2} x_{2n-3} x_{2n-4} x_{2n-5}]=e$\\
(b) $[x_{2n-7}, x_{2n-5}^{-1} x_{2n-4}^{-1} x_{2n-3}^{-1} x_{2n-2}^{-1} x_{2n-1}^{-1} x_{2n} x_{2n-1} x_{2n-2} x_{2n-3} x_{2n-4} x_{2n-5}]=e$\\
(c) $[x_{2n-6}, x_{2n-5}^{-1} x_{2n-4}^{-1} x_{2n-3}^{-1} x_{2n-2}^{-1} x_{2n-1} x_{2n-2} x_{2n-3} x_{2n-4} x_{2n-5}]=e$\\
(d) $[x_{2n-7}, x_{2n-5}^{-1} x_{2n-4}^{-1} x_{2n-3}^{-1} x_{2n-2}^{-1} x_{2n-1} x_{2n-2} x_{2n-3} x_{2n-4} x_{2n-5}]=e$
\end{tiny}

\medskip

We will show that these conjugations can be simplified.

By relation (3) of triple T4 and relation (3) of triple T3, relations (c) and (d) can be written as:\\
(c) $[x_{2n-6}, x_{2n-5}^{-1} x_{2n-1} x_{2n-5}]=e$\\
(d) $[x_{2n-7}, x_{2n-5}^{-1} x_{2n-4}^{-1} x_{2n-1} x_{2n-4} x_{2n-5}]=e$

\medskip

By relation (3) of triple T2 for $k=n-3$, we have $[x_{2n-6}, x_{2n-5}]=e$, and hence relation (c) becomes $[x_{2n-6}, x_{2n-1}]=e$.

By relation (1) of triple T2 for $k=n-3$, we have: $$x_{2n-7} x_{2n-5} x_{2n-4} =x_{2n-5} x_{2n-4} x_{2n-7} = x_{2n-4} x_{2n-7} x_{2n-5},$$ and then relation (d) becomes:
$[x_{2n-7}, x_{2n-1}]=e$.

\medskip

Now, we simplify relation (a). Again, by relations (2) and (3) of triple T2 for $k=n-3$, we have $[x_{2n-6}, x_{2n-5}]=e$ and  $[x_{2n-6}, x_{2n-4}]=e$, and hence:
$$[x_{2n-6}, x_{2n-5}^{-1} x_{2n-3}^{-1} x_{2n-2}^{-1} x_{2n-1}^{-1} x_{2n} x_{2n-1} x_{2n-2} x_{2n-3} x_{2n-5}]=e$$
By relation (3) of triple T4, we have:
$$[x_{2n-6}, x_{2n-5}^{-1} x_{2n-1}^{-1} x_{2n-3}^{-1} x_{2n-2}^{-1} x_{2n} x_{2n-2} x_{2n-3} x_{2n-1} x_{2n-5}]=e.$$
By relations (1) and (2) of triple T4, we have:
$$[x_{2n-6}, x_{2n-5}^{-1} x_{2n-1}^{-1} x_{2n} x_{2n-1} x_{2n-5}]=e.$$
By relations (1) of triple T3, we finally have: $[x_{2n-6}, x_{2n}]=e$.

\medskip

Now, we simplify relation (b). By relation (3) of triple T4, we have:
$$[x_{2n-7}, x_{2n-5}^{-1} x_{2n-4}^{-1} x_{2n-1}^{-1} x_{2n-3}^{-1} x_{2n-2}^{-1} x_{2n} x_{2n-2} x_{2n-3} x_{2n-1} x_{2n-4} x_{2n-5}]=e.$$
By relations (1) and (2) of triple T4, we get:
$$[x_{2n-7}, x_{2n-5}^{-1} x_{2n-4}^{-1} x_{2n-1}^{-1} x_{2n} x_{2n-1} x_{2n-4} x_{2n-5}]=e.$$
By relations (2) and (3) of triple T3, we get:
$$[x_{2n-7}, x_{2n-5}^{-1} x_{2n-4}^{-1} x_{2n} x_{2n-4} x_{2n-5}]=e.$$
Finally, by relation (1) of triple T2 for $k=n-3$, we have
$$x_{2n-7} x_{2n-5} x_{2n-4} =x_{2n-5} x_{2n-4} x_{2n-7} = x_{2n-4} x_{2n-7} x_{2n-5},$$
so we get: $[x_{2n-7}, x_{2n}]=e$.

\medskip

By similar tricks, one can simplify the conjugations for all the cases where $j=n$ and $i \leq n-4$. Hence, we have a presentation based on the topological generators without conjugations in the relations, and hence we are done.

\medskip

The above proof is based on the fact that the multiplicity of each multiple point is $3$. We now explain why it can be generalized to any multiplicity. In case of higher multiplicities, the quadruples from the previous case will be transformed to a block of $(n-1)(m-1)$ simple points. It can be easily checked that all the conjugations can be simplified in this case. Moreover, the triples from the previous case will be transformed into blocks similar to the blocks we had in the case of a cycle of length $3$ (see Proposition \ref{triangle-prop}), and in this case too, it can be easily checked that all the conjugations can be simplified, and hence we have shown that arrangements whose graph is a unique cycle and have no line with more than two multiple points, have a conjugation-free geometric presentation.
\hfill \qed

\bigskip

Using the following decomposition theorem of Oka and Sakamoto \cite{OkSa}, we can generalize the result from the case of one cycle to the case of a union of disjoint cycles:

\begin{thm} {\bf (Oka-Sakamoto)} Let $C_1$ and $C_2$ be algebraic plane curves in $\C ^2$.
Assume that the intersection $C_1 \cap C_2$
consists of distinct $d_1 \cdot  d_2$ points, where $d_i \ (i=1,2)$ are the
respective degrees of $C_1$ and $C_2$. Then:
$$\pi _1 (\C ^2 - (C_1 \cup C_2)) \cong \pi _1 (\C ^2 -C_1) \oplus \pi _1 (\C ^2 -C_2)$$
\end{thm}

Hence, we have the following result:
\begin{cor}
If the graph of the arrangement $\mathcal L$ is a union of disjoint cycles of any length and the arrangement has no line with more than two multiple points, then its fundamental group has a conjugation-free geometric presentation.
\end{cor}

\section*{Acknowledgments}
We would like to thank Patrick Dehornoy, Uzi Vishne and Eran Liberman for fruitful discussions.

We owe special thanks to an anonymous referee for many
useful corrections and advices and for pointing out the connection between our presentations and Dehornoy's positive presentations (Remark \ref{rem_dehornoy}).


\begin{thebibliography}{99}

\bibitem{AB} E. Artal-Bartolo, {\it Fundamental group of class of rational cuspidal curves}, Manuscripta Math. {\bf 93}, 273--281 (1997).

\bibitem{AB1} E. Artal-Bartolo, {\it A curve of degree five with non-abelian fundamental group}, Topology Appl. {\bf 83}, 13--29 (1997).

\bibitem{AB-CR} E. Artal-Bartolo and J. Carmona-Ruber, {\it Zariski pairs, fundamental groups and Alexander polynomials}, J. Math. Soc. Japan {\bf 50}(3), 521--543 (1998).

\bibitem{ABCT} E. Artal-Bartolo, J.I. Cogolludo and H. Tokunaga, {\it A survey on Zariski pairs},
in: {\it Algebraic geometry in east Asia, Hanoi 2005}, Adv. Stud. Pure Math. {\bf 50}, Math. Soc. Japan, Tokyo, 1--100 (2008).

\bibitem{arvola} W. Arvola, {\it The fundamental group of the complement of an arrangement of complex hyperplanes}, Topology {\bf 31}, 757--765 (1992).

\bibitem{chisini} O. Chisini, {\it Sulla identit{\`a} birazionale delle
    funzioni algebriche di due variabili dotate di una medesima curva
    di diramazione}, Rend. Ist. Lombardo {\bf 77}, 339--356 (1944).

\bibitem{CS} D.C. Cohen and A.I. Suciu, {\it The braid monodromy of plane algebraic curves and hyperplane arrangements},
    Comment. Math. Helv. {\bf 72}, 285--315 (1997).

\bibitem{CS-Hom} D.C. Cohen and A.I. Suciu, {\it Homology of iterated semidirect products of free groups}, J. Pure Appl. Alg. {\bf 126}, 87--120 (1998).

\bibitem{Deg} A.I. Degtyarev, {\it Quintics in $\C\PP^2$ with nonabelian fundamental group},
Algebra i Analiz {\bf 11}(5), 130--151 (1999) [Russian]; English translation:
St. Petersburg Math. J. {\bf 11}(5), 809--826 (2000).

\bibitem{Deh1} P. Dehornoy, {\it Groups with a complemented presentation}, J. Pure Appl. Algebra {\bf 116}, 115--137 (1997).

\bibitem{Deh2} P. Dehornoy, {\it Complete positive group presentations}, J. Algebra {\bf 268}, 156--197 (2003).

\bibitem{Deh3} P. Dehornoy, {\it The subword reversing method}, preprint (2009) [Available online: http://arxiv.org/abs/0912.4272].

\bibitem{ELST} M. Eliyahu, E. Liberman, M. Schaps and M. Teicher, {\it Characterization of line arrangements for which the fundamental group of the complement is a direct product}, Alg. Geom. Topo., to appear.  [Available online: http://arxiv.org/abs/0810.5533].


\bibitem{FN} E. Fadell and L. Neuwirth, {\it Configuration spaces}, Math. Scand. {\bf 10}, 111--118 (1962).

\bibitem{Fa1} K.M. Fan, {\it Position of singularities and fundamental group of the
   complement of a union of lines}, Proc. Amer. Math. Soc. {\bf 124}(11),
   3299--3303 (1996).

\bibitem{Fa2} K.M. Fan, {\it Direct product of free groups as the fundamental
   group of the complement of a union of lines}, Michigan Math. J.
   {\bf 44}(2), 283--291 (1997).

\bibitem{Fa3} K.M. Fan, {\it On parallel lines and free group}, unpublished note (2007) [online: http://arxiv.org/abs/0905.1178].

\bibitem{FrTe} M. Friedman and M. Teicher, {\it On non fundamental group equivalent surfaces}, Alg. Geom. Topo. {\bf 8}, 397--433 (2008).

\bibitem{Ga} D. Garber, {\it On the connection between affine and projective
fundamental groups of line arrangements and curves},
Singularit\'es Franco-Japonaises (J.-P. Brasselet and T. Suwa,
eds.), S\'eminaires \& Congr\`es {\bf 10}, 61--70 (2005).

\bibitem{GaTe} D. Garber and M. Teicher, {\it The fundamental group's
    structure of the complement of some configurations of real line
    arrangements}, Complex Analysis and Algebraic Geometry,
    edited by T.~Peternell and F.-O.~Schreyer, de Gruyter, 173-223 (2000).

\bibitem{GTV} D. Garber, M. Teicher and U. Vishne, {\it Classes of wiring diagrams and their invariants}, J. Knot Theory Ramifications {\bf 11}(8), 1165--1191 (2002).

\bibitem{Go} J.E. Goodman, {\it Proof of a conjecture of Burr, Gr\"unbaum and
   Sloane}, Discrete Math. {\bf 32}, 27--35 (1980).

\bibitem{HR} D.F. Holt and S.E. Rees, {\it The isomorphism problem for
   finitely presented groups}, in: {\it Groups, Combinatorics and Geometry},
   London Math. Soc. Lect. Notes Ser. {\bf 165},  459--475 (1992).

\bibitem{JY} T. Jiang and S.S.-T. Yau, {\it Diffeomorphic types of the complements of arrangements of hyperplanes}, Compositio Math. {\bf 92}(2), 133--155 (1994).

\bibitem{vK} E.R. van Kampen, {\em On the fundamental group
    of an algebraic curve}, Amer. J. Math. {\bf 55}, 255--260 (1933).


\bibitem{Kul} V.S. Kulikov, {\it On Chisini's conjecture},
    Izv. Ross. Akad. Nauk Ser. Mat. {\bf 63}(6), 83--116 (1999) [Russian];
    English translation: Izv. Math. {\bf 63}, 1139--1170 (1999).

\bibitem{Kul2} V.S. Kulikov, {\it On Chisini's conjecture II},
    Izv. Ross. Akad. Nauk Ser. Mat. {\bf 72}(5), 63--76 (2008) [Russian];
    English translation: Izv. Math. {\bf 72}(5), 901--913 (2008).

\bibitem{KuTe} V.S. Kulikov and M. Teicher, {\it Braid monodromy
    factorizations and diffeomorphism types},
    Izv. Ross. Akad. Nauk Ser. Mat. {\bf 64}(2), 89--120 (2000) [Russian];
    English translation: Izv. Math. {\bf 64}(2), 311--341 (2000).


\bibitem{milnor} J. Milnor, {\it Morse Theory}, Ann. Math. Stud.
    {\bf 51}, Princeton University Press, Princeton, NJ (1963).

\bibitem{BGT1}
 B. Moishezon and M. Teicher,  {\em Braid group technique in complex geometry I:
Line arrangements in $\C\P^2$}, Contemp. Math. {\bf 78}, 425--555 (1988).

\bibitem{BGT2}
B. Moishezon and M. Teicher,  {\em Braid group technique in complex geometry II: From arrangements of lines and
conics to cuspidal curves}, in Algebraic Geometry, Lect. Notes in Math. {\bf 1479},  131--180 (1990).

\bibitem{OkSa}
M. Oka and K. Sakamoto, {\it Product theorem of the
fundamental group of a reducible curve}, J. Math. Soc. Japan {\bf 30}(4), 599--602 (1978).

\bibitem{OT} P. Orlik and H. Terao, {\it Arrangements of hyperplanes}, Grundlehren der Mathematischen Wissenschaften {\bf 300},
    Springer-Verlag, Berlin, 1992.

\bibitem{Ra} R. Randell, {\it The fundamental group of the complement of a union of complex hyperplanes}, Invent. Math. {\bf 69}, 103--108 (1982).
    {\it Correction}, Invent. Math. {\bf 80}, 467--468 (1985).

\bibitem{Sal} M. Salvetti, {\it Topology of the complement of real hyperplanes in $\C^N$}, Invent. Math. {\bf 88}, 603--618 (1987).

\bibitem{WY} S. Wang and S.S.-T. Yau, {\it Rigidity of differentiable structure for new class of line arrangements}, Comm. Anal. Geom. {\bf 13}(5), 1057--1075 (2005).

\bibitem{Z1} O. Zariski, {\it On the problem of existence of
    algebraic functions of two variables possessing a given
    branch curve}, Amer. J. Math. {\bf 51}, 305--328 (1929).

\bibitem{Z2} O. Zariski, {\it On the Poincar\'e group of rational plane curves},
   Amer. J. Math. {\bf 58}, 607--619 (1936).


\end{thebibliography}
\end{document}